# Kazakh History and Philosophy: the Ethnomathematical Component of the Content of Primary School Education in the Republic of Kazakhstan

**Nurassyl Kerimbayev, Aliya Akramova**
Almaty, Kazakhstan
N_nurassyl@mail.ru

**Abstract**
In the article, we would like to take part in a discussion over the problem of ethnomathematics, which has been under way for several years. We would like to express our opinion, our ideas concerning this problem and suggest our point of view on solving this problem. We would like to introduce readers to experience in working on including ethnomathematical component when teaching mathematics at primary school in the Republic of Kazakhstan.

## I. Important points

The most important question, which scientists ask, is the following: can mathematics be "ethnic"? A lot of research consider ethnomathematics in connection with the variety of cultures and human everyday activities. Ethnomathematics is associated with human factor.

The Republic of Kazakhstan is a tolerant, multicultural state (there are over 130 nations living on its territory). We have raised the question: how can equilibrium between social and political aspects of math education be maintained in these conditions? How can we save the originality of ethnic groups, absorb values and standards of different cultures when acquiring mathematics as a branch of science?

It is necessary to note, that in the Republic of Kazakhstan at the top government executive level the "Concept of ethnic cultural education" has been adopted. The main idea of the Concept is "realizing ethnic cultural interests of inhabitants in the sphere of mathematics". That very idea provides a legal basis for a democratic society.

The importance of our research is that at math lessons problems of acquiring a necessary level of math education and social and cultural problems were solved.

## Introduction

In multi-cultural regions with pupils of different nations and languages, there arises necessity to take into account ethnocultural and ethnic component. In our article, we consider the ethnomathematical component in the content of mathematical education ***not as*** ethnic or mathematical knowledge of some territorial area. We consider the ethnic component with relation to philosophical and historical aspects. The Kazakh history and philosophy contain a copious ethnic material, which has a cognitive, pedagogical and educational potential. That is why we consider historical and mathematical material to be the ethnomathematical component in the content of math education at lessons in primary school.

Ethnomathematics and the definition of this term are the object of many scientists' research. Ubiratan D'Ambrosio (2002) defines ethnomathematics as a "scientific research program in the sphere of history and philosophy of mathematics". He states that ethnomathematics has "evident pedagogical consequences".

Problems of ethnomathematical education and pedagogical knowledge, including roles of teachers in this process, have been discussed in the articles by Oslund, Joy A. (2012)*,* Pais, Alexandre, and Paola Valero (2014)*,* de Freitas, Elizabeth, and Nathalie Sinclair (2013) Brown, Tony, and Margaret Walshaw (2012), they position math education with relation to a political practice. In the article «Researching research: mathematics education in the Political» they revealed Fuko's concept of bio-politics with respect to the object "Mathematics" (Pais, Alexandre, and Paola Valero, *2012*). Stinson, David W., and Erika C. Bullock provide an overview of mathematics education as a research domain, identifying and briefly discussing four shifts or historical moments (2012).

*Bill Barton* emphasizes complexity of defining the term and difficulties arising when studying this problem: «There is little agreement on the extent to which mathematics is universal, and on how mathematical ideas can transcend cultures. Very little of the ethnomathematical literature is explicit about its philosophical stance. This is one of the areas which must be addressed if the subject is to gain wider legitimacy in mathematical circles» *(1996)*. Vithal, Renuka, and Ole Skovsmose explores a critique of ethnomathematics using the South African situation and conceptual tools of a critical mathematrics education (1997). Rowlands, Stuart, and Robert Carson give response to Adam, Alangui and Barton's ``A Comment on Rowlands & Carson `Where would Formal, Academic Mathematics stand in a Curriculum informed by Ethnomathematics? A Critical Review'' (2004). Knijnik, Gelsa discusses a new philosophical perspective for ethnomathematics which articulates Ludwig Wittgenstein's *(2012)*. Vilela, Denise Silva considers a philosophical background for the ethnomathematical program (2010). The object of study of Pinxten, Rik, and Karen François *(2011) is* bicultural school in Navaho reservation in Arizona. De Freitas Elizabeth, and Nathalie Sinclair (2013) put forward new materialist ontologies in mathematics education. The authors consider math concepts to be a dynamic material.

Problems of ethnic minorities, demographic, linguistic features of different groups and sections of the population in education are widely discussed. Stein, Robert George (2013) introduces to the history of school mathematics in colonial North America in 1607-1861. The Greek national-linguistic community of Italy, increase of exogamy have been studied by Biondi, G., and E. Perrotti (1990). *The results* of the educational program of 14 ethnic groups Shingu of the Indian park (Brazil). In the research they consider the textbook on Mathematics written in the language of the native population, which will be used in the schools of the park (Mendes, Jackeline Rodrigues, 2007).

The article by Roth, Wolff-Michael tells of the importance of a personality, its subjective factors in math education. The article considers principle ideas, which L.S. Vygotskiy put forward. Theory of personality and activity is basic for most investigations in the sphere of teaching and educating. *(Roth, Wolff-Michael, 2012)*

In our work we analyzed the textbooks on Mathematics for the first – forth forms of secondary schools. The article «Opening mathematics texts: resisting the seduction» by Wagner, David helped us to construct our analysis. The author's basic idea is importance of modeling study material and critical attitude to mathematical texts. (Wagner, David, *2012)*.

Adam Shehenaz, Wilfredo Alangui, and Bill Barton in the article (2003) state that "The critique of ethnomathematics that appeared recently provides an opportunity to open debate on cultural issues in mathematics. This response argues that such debate must be based on contemporary writing in the field, and should not focus on extreme views within the political justification for ethnomathematics. It addresses some of the philosophical questions raised by Rowlands and Carson, and the relationship of the field with indigenous knowledge is raised. We also suggest that the role of ethnomathematics in mathematics education is now predominantly an empirical matter, and comment on some preliminary results from recent studies that indicate **a positive role** for culturally-based curricula".

Our experience and the results of the research carried out showed a positive effect of the ethnic component on the content of math education.

**II. Methodology of research**

The main methodological tools were a theoretical analysis of the literature survey and interviews, the analysis of the results, forecasting the future direction of research.

In order to include elements of ethnomathematics we have studied folkloristics, the literature on philosophy of the Kazakh people. We got acquainted with the mathematical material, which is related to such subjects as quantities, arithmetical and geometrical material. This information was given as brief historical references. The fact is that many mathematical facts relate to the life, culture and philosophy of the Kazakh people, and these facts are unknown for a wide range of other people.

We carried out questionnaire, testing among primary school students, teachers and parents. We published our testing results in our scientific articles in some journals (2014a, 2014b, 2014c). We wanted to learn social opinion concerning including ethnomaterial in the process of teaching mathematics at primary school. Including ethnomaterial in Mathematics textbooks for 2 – 4 forms was welcomed by teachers, because this information had been unknown or forgotten. The statistical analysis of the experiment results and their methodical interpretation demonstrated a positive change in the theory and methods of teaching mathematics at primary school.

### III. Body
### III.I What do we mean by ethnomathematical component in the content of math education?

Present researchers associate ethnomathematics with the cultural environment of an ethnos where math ideas are described using the practice of these cultures. The main idea of our including an ethnomathematical material is acquiring the culture of the Kazakh people. Using elements of the Kazakh folklore and Kazakh national philosophy is a condition of getting to know the culture of the Kazakh people and integrating in other cultures. So, we can speak of moving from a small "micro-objective" of the ethnocultural component to a more global one, which fulfills a more important mission in social and cultural relations of various groups and sections of the population (Figure 1).

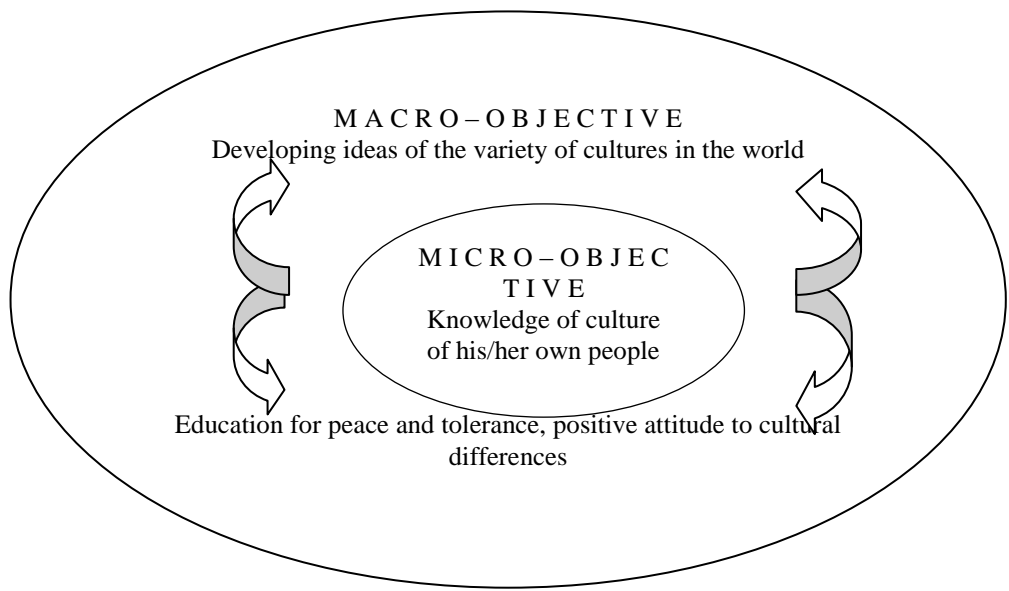

Figure 1 – Linkage between micro- and macro-objectives when developing an ethnomathematical component in the content of math education

In the article we suggest getting acquainted with some ideology of the Kazakh people. The philosophic view of the Kazakh people is included in the processes of world experience and formed the basis of many mathematical ideas.

Philosophy of the Kazakh people has developed under the influence of centuries-old observation of phenomena of the world around, changes in nature, social and historical events. The centuries-old experience of the people was reflected in legends, tales, proverbs and sayings. There is the whole experience of mankind cognition in the Kazakh folklore. For the Kazakh society philosophizing in non-philosophic forms at the level of universals of ideology, which forms human activities, is typical. For the Kazakhs Nature has always been an absolute and eternal beginning, the only reality, which has immanent reasons and factors of existence, basics of activity and evolution. At first people could only obey natural laws, understand their reasons and develop themselves using them. Here nature performs as the Universe, having inexhaustible depths and logic of its existence, which dictates to individuals forms and ways of their existence.

For instant, when measuring time, in many languages they say: morning, afternoon, evening, night. Besides, we can say: sunrise, early morning, late morning (anyway, morning); early afternoon, late afternoon (anyway, afternoon); early evening, late evening (anyway, evening). There is also just after noon, noon, and twilights.

How do the Kazakhs differ day and nighttime? Let us consider some periods of day and nighttime the Kazakhs use. Each period of time is divided into 3 – 5 time periods and has its own name (Table 1):

Table 1.

| I. **Tan erten** — morning | II. **Saske** — the time when the Sun rises above the horizon | III. **Tus** — afternoon | IV. **Besin** — after noon, the Sun crossed the zenith | V. **Aksham, ymyrt** — evening, twilights | VI. **Tun** — night |
|---|---|---|---|---|---|
| а) elen – alan — at dawn, day break | а) siyr saske — the time, when the Sun has risen at the length of the lasso | а) tal tus, talma tus, tapa-tal tus — approximately at 1 o'clock in the afternoon, day peak | а) uly besin — the Sun begins moving far down the west | а) alageuim — the Sun is on its way to set, early evening dusk | а) inir — the time before nightfall, the nature is getting ready to sleep |
| б) kulaniek, kulancary — when down breaks and one can make out the shapes of objects | б) sasketus — approximately at 12 o'clock in the afternoon | б) shankai tus — the daytime when shadow is the shortest, approximately at 2 o'clock in the afternoon | б) kishi besin — the Sun sank visibly | б) keugim, keuim, ymyrt — the Sun set, dusk deepens | б) kyzyl inir — the beginning of night |
| в) tansary — the time of grey light, the Sun has not risen yet | в) uly saske — noon, midday | | в) kulama besyn — the Sun sank more | в) keshkurym, namazsham — the time of evening prayers | в) zharym tun — midnight |
| г) tan — the time when the Sun appears, rises | | | г) yekindi — the evening is coming, the Sun is low | г) kesh — everything is steeped by twilight, the beginning of night | г) tan karangysy — darkness, the time just before sunrise |
| | | | д) namazdyger — the Sun is getting ready to set | | |

This material was offered to children as an introductory historical material on pages of "Mathematics" textbook when they studied time units.

> **Time Units**
> **The old Kazakh time units were the following:** sut pisirim, biye sauym, yet pisirim, saske tus, shankai tus.
> Sut pisirim – the time for boiling milk (5–10 minutes).
> Biye sauym – the perifod of time between milking a mare (1–2 hours).
> Yet pisirim – the time for cooking meat (3–4 hours).
> Some periods of daytime were called saske tus, shankai tus.
> Saske tus – the time of a day, when there is no shadow (12 o'clock).
> Shankai tus – the time of a day, when shadow is the shortest (2 o'clock in the afternoon).

Nomadism necessitated estimating distances, traveled way. The measures of length of the Kazakh people also have a philosophic character, which reflects nomads' natural watchfulness.

Ontological characteristics of quantity have always been space and time. Proportions were separated in the process of measuring, and this was used when counting number of domestic animals (cattle) or for land allotment. The Kazakhs' labour activity resulted in developing the conception of space as the distance of roaming from one place to another within 100 km and more, which they passed for a certain number of days and months, and they used distances of day roaming (kosh), week roaming, etc.

On wintering areas and settlements there were delineated land plots separated by aryks, which had from 60 to 150 square meters, depending on numbers of members of families. It demonstrates geometrical point of view at space and functions of measuring it. It reflected optimal acres in crops on irrigated land. The concept of quantity was developed in the systems of cattle-rearing and farming in the same way, which is evidence of a single symbolic system of the Kazakhs, in which the code of the social experience of the people is registered and transferred. Beginning with the counting operation there is development of basic sub-categories of quantity – numbers, magnitudes, structures.

**Mass measures**
**Kazakh old mass units were** pitir bidai, zhamby kumis, korzhyn, arka zhuk, dagar.
1 pitir bidai – 3 kilos of wheat.
1 zhamby kumis – silver piece weighing 6 kilos.
Korzhyn – the weight of a bag with 2 sections (40 kilos).
1 arka zhuk – the weight of the luggage (45–50 kilos), which one can carry on the back.
Dagar – the weight of a bag with seeds (100 kilos).
**Capacity measures**
**Kazakh old measures of capacity were** shara, torsyk, konek, saba.
Shara – a bowl for koumiss (2 liters).
Torsyk – a leather dish (6–8 liters).
Konek – a pail for milking a mare (15–16 liters).
Saba – dish for koumiss (10–15 liters).

When forming an ethnocultural component in the process of teaching Mathematics at primary school the approach to education is realized not just as the mechanism of transferring knowledge but also as a culture-generating component, the most important means of saving and developing human and national originality by an individual (Figure 2).

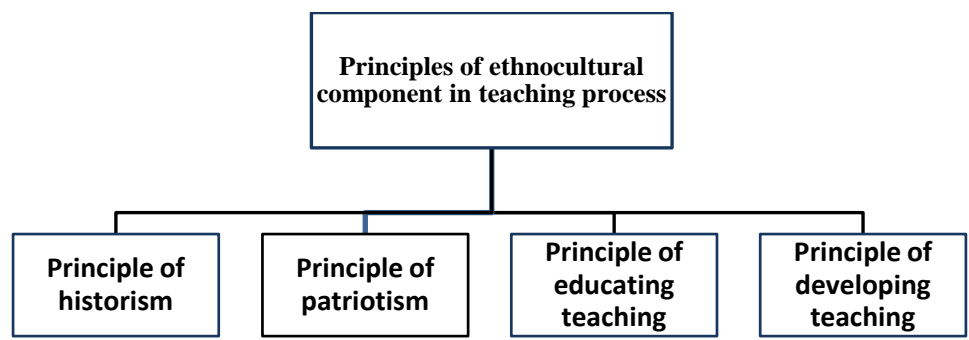

Figure 2. Principles of ethnocultural component in the process of teaching Mathematics.

We mark out the following ways of developing and realizing the ethnic component in the process of teaching Mathematics to primary school students (Figure 3):

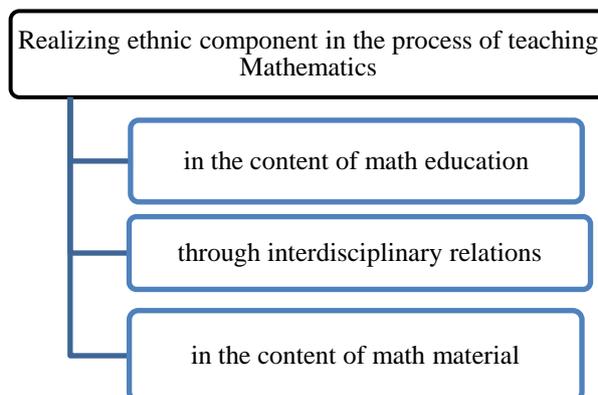

Figure 3. Realizing ethnic component in the process of teaching Mathematics.

Summarizing mentioned above, we can conclude that psychological and pedagogical basics of ethnocultural work at primary school Mathematics lessons propose the following:

- Taking into account specific features of primary school students; activating students' motivation for curricular and extra-curricular activity;
- Improving pedagogical process; integrating teaching and educating;
- Creating an educational environment on the basis of reviving traditions, customs, rites of our people, elements of national pedagogy; создание обучающей среды основе возрождения народных традиций, обычаев, обрядов, элементов народной педагогики;
- Increasing teachers' competence and professionalism, meaning acquiring new pedagogical technologies based on Kazakh folklore and Kazakh philosophy.

### III.II Analysis of the Program and educational-methodical complexes on Mathematics for primary school in ethnoculatural aspect

Realization of the principle of ethnocultural education based on Kazakh folklore and Kazakh philosophy in the process of math education proposes a motivated, systematic use of the ethnocultural component, the structural elements of which are ethnic culturological study material selected in accordance with the content of the State Compulsory Standard of Education of the RK and Program on Mathematics for Primary school, where the Kazakh folklore is reflected (Kucher T.P, Akramova A.S., Kukarina G.I., 2013).

The analysis of the Program on Mathematics on developing ethnocultural competence allowed us to identify the ethnocultural aspect of the program material for 1, 2, 3, 4 forms (Table 2):

*Table 2.* Analysis of the Program on Mathematics on developing ethnocultural competence.

| Form | Ethnocultural aspect |
|---|---|
| 1 form | Cultivating the culture of communication, respect for the older and care for the younger; active participation in environmental protection; patriotism; respect for the history, culture, traditions and other values of the Kazakh people and other nations living in Kazakhstan. |
| 2 form | Brief historical information is provided for familiarizing with achievements of the domestic and international culture as an introductory material. For instance, how ancient people used to write equations, how different ancient nations used to write single-digit, two-digit, three-digit numbers. Why is the number 40 called "forty" but not "four ten"? What is Al-Farabi? What is Pythagoras? The history of the magic square. Ancient Kazakh units of mass, units of length, units of capacity, unites of time. |
| 3 form | Developing the idea of mathematics as a part of general human culture, understanding the importance of mathematics for social progress; Cultivating respect for human rights and liberties; culture of communication, respect for the older and care for the younger; culture of health, ability of social adaptation; civism, love for the surrounding nature, homeland, and family. |

| 4 form | For familiarizing students with achievements of domestic and international culture, at some lessons some brief historical information on Mathematics related to the content of the material under study is given as an introductory material. For instance, how different ancient people used to write multidigit numbers. How and when did the modern number writing appear? The German scientist Karl Gauss. What is percent? Parallel straight lines. From the history of computing devices. Which of the old units do we use today? When did fractions appear? The history of appearance of month names. From the history of clocks and watches. |

In the Program on Mathematics (elementary school, 1–4 forms) development of the ethnocultural component is considered as a personal result of study (Table 3):

*Table 3.* Analysis of the Program requirements for the level of students' knowledge in ethnocultural aspect.

| Forms | Ethnocultural aspect |
|---|---|
| 1–4 forms | **3.1.2 Personal results:** Respect for the history, culture and traditions and other values of the Kazakh people and other nations living in Kazakhstan. |

### III.III Methodological particularities of including an ethnocultural component in the process of teaching Mathematics to elementary school students on the basis of Kazakh philosophy

The ethnocultural component in the process of teaching Mathematics to elementary school students leads to the possibility to consider phenomena in their historical retrospect; helps reveal relations of the phenomena under study on the basis of Kazakh philosophy.

The culture of each ethnos is peculiar, original, and unique. When developing an ethnocultural component at Mathematics lessons it is necessary to use the method of empathy, sympathy, diving into the logic of other people.

A necessary condition for developing an ethnocultural component in the content of math education is **creation** of ethnocultural educational space (http://adilet.zan.kz/rus/docs/N960003058).

The love for one's people should go together with the love for the land where the people lives. The basic part of ethnocultural education is institutional education. An ethnocultural educational space is created by the following aspects: in which language subjects are taught, which subjects are taught, and what content of the courses under study is.   The result of such education is developing a multi-cultural person. A multi-cultural person means an individual oriented at other cultures through his own culture. A deep knowledge of his own culture for him is a basis of interest in others cultures, and familiarizing with many cultures is a basis of spiritual enrichment and development.

Multi-cultural education has a positive effect on the social and psychological environment, in which every student independently on his/her identity has the same possibilities as others for realizing his/her constitutional rights for getting equal education, for realizing his/her potential and social development while studying. Such education has extensive possibilities for reducing conflicts in a community, developing tolerance.

An ethnocultural component of math education can be presented in the form of various ethnic material (Figure 4).

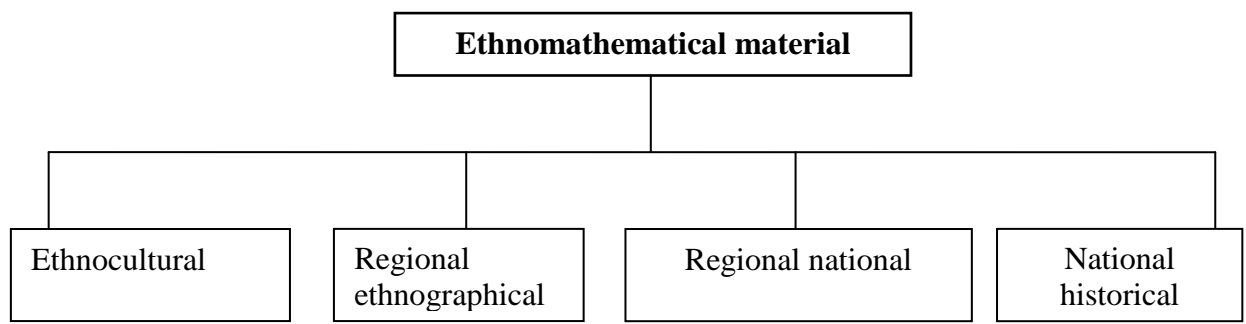

*Figure 4.* Types of ethnomathematical material

Solving different math problems of ethnic character provides with developing dialectic-mathematical understanding of nature, broadens the mind, connects Mathematics with surrounding reality through Kazakh philosophy.

**IV. Results**

The experimental work was aimed at testing the effectiveness of the methods of developing elementary school students' ethnocultural competence at Mathematics lessons.

For conducting the experiment the parallel test groups: experimental and control (according to equaling the groups as a whole) were chosen. The original level of elementary school students' developed ethnocultural competences was defined by a number of diagnostic tasks. At the end of the experiment, we studied a summary level of development of the competences.

During the experimental teaching, we paid a special attention to including the tasks aimed at developing ethnocultural competence of elementary school students at Mathematics lessons on the basis of Kazakh folklore and Kazakh philosophy.

The teaching experiment was conducted on the two lines:
1. Teaching children in the control class was conducted traditionally.
2. Teaching children in the experimental class was conducted on the basis of the system of tasks developed, which contained an ethnocultural component and provided with developing an ethnocultural competence.

The experiment was realized in 3 stages. At the first stage, the focus of the experimental teaching was determining the level of existing ethnocultural competences using special tasks.

At the second stage a special attention was paid to applying a system of problems on Mathematics aimed at developing ethnocultural competence. At the third stage of the experimental teaching, we continued to improve solving problems on developing ethnocultural competence.

We carried out our experiment in schools № 51 and № 73 of Almaty. At the first stage we gave testing, diagnostic tasks, talked to the children. All the problems were based on the same math concepts. The difference was only in the conditions of applying them, in novelty of the pedagogical problems of the research, in the didactical material, which served as visual aids.

The control material was used during the experiment for determining the level of the children's gaining the material at the beginning of the experiment and after each stage, for considering laws of acquiring the material by a child, for analyzing the conditions providing with the most effective teaching elementary school students.

During our experiment the following principles we drew on played an important role:
- ethnocultural level of education;
- succession;
- personal creative potential;
- taking into account students' ethnopsychological and psychophysiological features

As criteria, we determined cognitive, emotional and activity-related criteria of development of children's ethnic competence at Math lessons.

The cognitive criterion is the result of a cognitive activity at Mathematics lessons including the knowledge of nationalities of the human society, variety of national games, folklore of different peoples. The emotional criterion is an interest in gaining knowledge of other nationalities and their cultures. The activity-related criterion is an active participation in solving math problems and carrying out ethnocultural exercises.

For each component we defined the levels of development of ethnic competence: high, average and low (Table 4).

*Table 4.* Levels of development of the components of ethnic competence.

| Levels | Characteristics |
|---|---|
| III level<br>High level | Children's answers are correct and complete, represent the facts, which demonstrate the knowledge of the culture of their people and cultures of other ethnic groups |

| II level Average level | Children's answers do not always represent the facts, are incomplete, are not argued; the knowledge of their culture and cultures of other ethnic groups is not complete |
| I level Low level | Children's answers do not often represent the facts, are brief, are not argued; there is very poor knowledge of their culture and cultures of other ethnic groups |

The first stage of the experiment carried out among the second form students from school № 51 demonstrated the following levels of development of cognitive, emotional and activity-related components (Table 5):

*Table 5.* Levels of development of ethnic competence of control and experimental classes students at the beginning of the experiment.

| CRITERIA | LEVELS | High | | Average | | Low | |
|---|---|---|---|---|---|---|---|
| | | Control class | Experimental cl. | Control class | Experimental cl. | Control class | Experimental cl. |
| | Cognitive | 40 | 45 | 35 | 35 | 25 | 20 |
| | Emotional | 45 | 50 | 35 | 40 | 20 | 10 |
| | Activity-related | 35 | 40 | 40 | 40 | 25 | 20 |

Thus, we can see from Table 4 that the high level of all the components prevails both in the control and experimental classes (cognitive is 45% (40%), emotional is 50 % (40%), and activity-related is 40% (35%)), which shows good results (Figure 5).

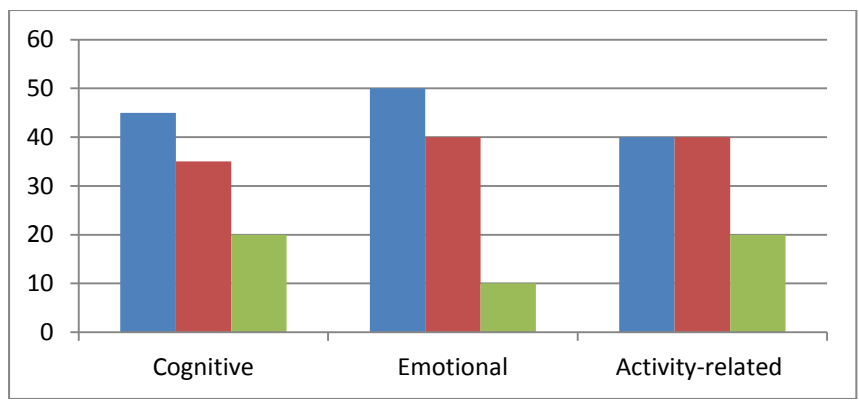

*Figure 5.* Levels of development of ethnic competence of control and experimental classes students at the beginning of the experiment.

The experimental data showed that most of the children have high and average levels concerning the level of development of basic components. During developing experiment (the second stage) we used various teaching aids, approbated a complex of exercises and tasks for developing ethnic competence of elementary school students at Mathematics lessons.

Use of the tasks on Mathematics containing an ethnocultural component improved students' knowledge of the life of their people, and other peoples and nations, elements of the Kazakh folklore broadened the mind of elementary school students. The third stage of the experiment demonstrated the following levels of development of the second-form students' cognitive, emotional and activity-related components (*Table 6,* Figure 6):

*Table 6.* Levels of development of ethnocultural competence at Mathematics lessons of control and experimental classes students at the end of the experiment.

| CRIT | LEVELS | High | | Average | | Low | |
|---|---|---|---|---|---|---|---|
| | | Control | Experimental | Control | Experimental | Control | Experimental |
| | Cognitive | 45 | 75 | 40 | 20 | 15 | 5 |
| E R I A | Emotional | 55 | 80 | 35 | 15 | 10 | 5 |
| | Activity-related | 40 | 70 | 40 | 20 | 20 | 20 |

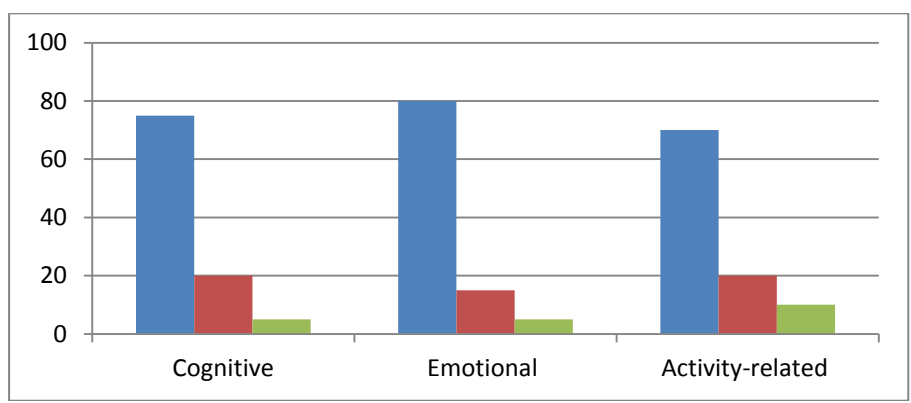

*Figure 6.* Levels of development of ethnocultural competence at Mathematics lessons of control and experimental classes students at the end of the experiment.

Therefore, we can state that the lessons on developing schoolchildren's ethnic tolerance in a multi-cultural space of elementary school, which we had selected, gave good results (Table 7).

*Table 7.* Time history of levels of development of ethnocultural competence of control and experimental classes students from schools № 51 and № 73 of Almaty by the end of the experiment

| Forms | | Levels of development of elementary school students' ethnic competence at Mathematics lessons | | | | | |
|---|---|---|---|---|---|---|---|
| | | High (%) | | Average (%) | | Low (%) | |
| | | Beginning | End | Beginning | End | Beginning | End |
| 2 form, schools № 51 and № 173 | Control | 40 | 35 | 40 | 40 | 20 | 15 |
| | Experimental | 45 | 75 | 35 | 20 | 20 | 5 |
| 3 form, schools № 51 and № 173 | Control | 50 | 55 | 40 | 40 | 20 | 5 |
| | Experimental | 50 | 65 | 35 | 30 | 15 | 5 |

.

According to the formula of the research efficiency coefficient, which is $R_{effect} = R_{experim.}/R_{control}$, we have: $R_{exper.}/R_{control} > 1$, which proves effectiveness of the chosen methods.

Thus, the hypothesis that use of special tasks, realization of the complete methodological system of work at Mathematics lessons at primary school contribute to develop ethnocultural competence on the basis of Kazakh philosophy and Kazakh folklore.

**V. Conclusion**

Development of ethnic component when teaching Mathematics is one of priority objectives of the educational process at primary school.

The content of an ethnomathematical material is based on the elements of Kazakh folklore and Kazakh philosophy. Such content can be widely used at Mathematics lessons when doing different exercises and solving different problems.

The research was carried out at theoretical and experimental levels. On the basis of the obtained data we determined a set of necessary pedagogical and methodical conditions of using a system of math problems containing elements of Kazakh folklore and Kazakh philosophy. These conditions provide with developing an ethnomathematical component at primary school Mathematics lessons.

The experiment we conducted demonstrated that elements and information of the Kazakh history and philosophy can be used as an ethnomathematical component of the content of primary education in the Republic Kazakhstan.